\documentclass[12pt]{article}
\usepackage{amsmath, amssymb, amsfonts, amsthm, amscd, vmargin}
\usepackage[latin1]{inputenc}
\usepackage[all]{xy}
\usepackage{graphicx}

\setmarginsrb{2.4cm}{4.5cm}{2.4cm}{3.3cm}{0cm}{0cm}{0cm}{1.5cm}

\theoremstyle{plain}
\newtheorem{teo}{Theorem}[section]
\newtheorem{lem}[teo]{Lemma}
\newtheorem{prop}[teo]{Proposition}

\theoremstyle{definition}
\newtheorem{defi}[teo]{Definition}

\newtheorem{rem}[teo]{Remark}
\newtheorem{rems}[teo]{Remarks}
\numberwithin{equation}{teo}
\newcommand{\Aut}{\operatorname{Aut}}

\newcommand{\Lie}{\operatorname{Lie}}

\newcommand{\Cbb}{{\mathbb C}}

\newcommand{\Zbb}{{\mathbb Z}}
\newcommand{\Rbb}{{\mathbb R}}
\newcommand{\Pbb}{{\mathbb P}}

\newcommand{\lra}{\longrightarrow}
\newcommand{\lmt}{\longmapsto}

\begin{document}

\title{Smooth projective horospherical varieties with Picard number~1}

\author{Boris Pasquier}
\maketitle

\begin{abstract} We describe smooth projective horospherical varieties with Picard number~1. Moreover we prove that the automorphism group of any such variety acts with at most two orbits and we give a geometric characterisation of non-homogeneous ones.\end{abstract}



\section*{Introduction}
A horospherical variety is a normal algebraic variety where a reductive algebraic group acts with an open orbit which is a torus bundle over a flag variety.  For example, toric varieties and flag varieties are horospherical.

In this paper we describe all smooth projective horospherical varieties with Picard number~1. In particular, they are Fano varieties,  {\it i.e.} their anticanonical bundle is ample. Fano horospherical varieties are classified in terms of some rational polytopes in \cite{Pa06}. 
The study of those with Picard number~1 is partly motivated by Theorem 0.1 of \cite{Pa06}, that gives an upper bound of the degree of smooth Fano horospherical varieties. Indeed, the expression of this upper bound is different if the Picard number is $1$. 

One can remark that the smooth projective toric varieties with Picard number~1 are exactly the projective spaces. This is not the case for horospherical varietes: the first examples of smooth projective horospherical varieties with Picard number~1 are flag varieties $G/P$ with $P$ a maximal parabolic subgroup of $G$. 

Moreover, smooth projective horospherical varieties with Picard number~1 are not necessarily homogeneous. For example, let $\omega$ be a skew-form of maximal rank on $\Cbb^{2m+1}$. For  $i\in\{1,\ldots,m\}$, define the odd symplectic grassmannian $Gr_{i,2m+1}$ as the variety of $i$-dimensional $\omega$-isotropic subspaces of $\Cbb^{2m+1}$. Odd symplectic grassmannians are horospherical varieties (see Proposition \ref{oddgrass}) and have two orbits under the action of their automorphism group which is a connected non-reductive linear group \cite{Mi05}.
 In fact, the latter property is satisfied by all non-homogeneous smooth projective horospherical varieties with Picard number~1 (see Theorem \ref{mainteo}).

Two-orbits varieties ({\it i.e.} normal varieties with two orbits under the action of a connected linear group) have already be studied by D.~Akhiezer and S.~Cupit-Foutou. Indeed, D.~Akhiezer has classified the ones whose closed orbit is of codimension~1 and he has proved in particular that  they are horospherical when the group is not semi-simple.  S.~Cupit-Foutou has classified two-orbits varieties when the acting group is semi-simple, and she has also proved that they are spherical ({\it i.e.} they admit a dense orbit of a Borel subgroup) \cite{CF03}. It could be interesting to look at smooth two-orbits varieties with Picard number~1.  Note also that a classification of smooth projective symmetric varieties with Picard number~1 has been recently given by A~Ruzzi \cite{Ru07}.\\

The paper is organized as follows. In Section \ref{1}, we recall the known results on horospherical varieties which we will use in all the paper. In particular we briefly summarize the Luna-Vust theory \cite{LV83} in the case of horospherical varieties.

In Section \ref{2}, we prove that a given horospherical homogeneous space admits at most one smooth compactification with Picard number~1. Then we give the list of horospherical homogeneous spaces that admit a smooth compactification not isomorphic to a projective space and with Picard number~1. We obtain a list of 8 cases (Theorem \ref{classif}).

In Section \ref{3}, we prove that in 3 of these cases, the smooth compactification is homogeneous (under the action of a larger group).

In Section \ref{4}, we study the 5 remaining cases. We compute the automorphism group of the corresponding smooth compactification with Picard number~1. We prove that this variety has two orbits under the action of its automorphism group and the latter is connected and not reductive.

In Section \ref{5}, we obtain a purely geometric characterization of the varieties obtained in Section \ref{4}.\\

\section{Notation}\label{1}

Let $G$ be a reductive and connected algebraic group over $\Cbb$, $B$ a Borel subgroup of $G$, $T$ a maximal torus of $B$ and $U$ the unipotent radical of $B$. Denote by $C$ the center of $G$ and by $G'$ the semi-simple part of $G$ (so that $G=C.G'$). Denote by 
$S$ the set of simple roots of $(G,B,T)$, and by $\Lambda$ (respectively $\Lambda^+$) the group of characters of $B$ (respectively the set of dominant characters).
Denote by $W$ the Weyl group of $(G,T)$ and, when $I\subset S$,  denote by $W_I$ the subgroup of $W$ generated by the reflections associated to the simple roots of $I$.
If $\alpha$ is a simple root, we denote by $\check\alpha$ its coroot and $\omega_\alpha$ the  fundamental weight corresponding to $\alpha$. Denote by $P(\omega_\alpha)$ the maximal parabolic subgroup containing $B$ such that $\omega_\alpha$ extends to $P(\omega_\alpha)$. Let $\Gamma$ be the Dynkin diagram of $G$. When $I\subset S$, we denote by $\Gamma_I$ the subgraph of $\Gamma$ with vertices the elements of $I$ and with edges those joining two elements of $I$.
Let $\lambda\in\Lambda^+$, then we denote by $V(\lambda)$ the irreducible representation of $G$ of heighest weight $\lambda$ and by $v_\lambda$ a heighest weight vector of $V(\lambda)$.

Let $H$ be a closed subgroup of $G$. Then $H$ is said to be {\it horospherical} if it contains the unipotent radical of a Borel subroup of $G$. In that case we also say that the homogeneous space $G/H$ is {\it horospherical}. Up to conjugation, one can assume $H\supset U$.  Denote by $P$ the normalizer $N_G(H)$ of $H$ in $G$, it is a parabolic subgroup of $G$ such that $P/H$ is a torus. Then $G/H$ is a torus bundle over the flag variety $G/P$. The dimension of the torus is called the {\it rank} of $G/H$ and denoted by $n$.

All varieties are irreducible algebraic varieties over $\Cbb$.
Let $X$ be a normal variety where $G$ acts. Then $X$ is said to be a {\it horospherical variety} if $G$ has an open orbit isomorphic to an horospherical homogeneous space $G/H$. In that case, $X$ is also said to be a $G/H$-embedding. The classification of $G/H$-embeddings (due in a more general situation\footnote{When $G/H$ is spherical. See for example \cite{Kn91} or \cite{Br97b}.} to D.~Luna et Th.~Vust \cite{LV83}) is detailed in \cite[Chap.1]{Pa06}. 

Let me summarize here the principal points of this theory.  Let $G/H$ be a fixed horospherical homogeneous space of rank $n$. Then it defines a set  of simple roots $$I:=\{\alpha\in S\mid \omega_\alpha \mbox{ is not a character of }P\}$$ where $P$ is the unique parabolic subgroup defined above.  We also introduce a lattice $M$ of rank $n$ as the sublattice of $\Lambda$ consisting of all characters $\chi$ of $P$ such that the restriction of $\chi$ to $H$ is trivial. Denote by $N$ the dual lattice to $M$.

For all $\alpha\in S\backslash I$, we denote by $\check\alpha_M$ the element of $N$ defined as the restriction to $M$ of $\check\alpha:\Lambda\lra\Zbb$. The point $\check\alpha_M$ is called the {\it image of the color} corresponding to $\alpha$. See \cite[Chap.1]{Pa06} for the definition of a color and note that the set of colors is in bijection with $S\backslash I$.

\begin{defi}
A {\it colored cone} of $N_\Rbb:=N\otimes_\Zbb\Rbb$ is an ordered pair $(\mathcal{C}, \mathcal{F})$ where $\mathcal{C}$ is a convex cone of $N_\Rbb$ and $\mathcal{F}$ is a set of colors called the set of colors of the colored cone, such that\\
(i) $\mathcal{C}$ is generated by finitely many elements of $N$ and contains the image of the colors of $\mathcal{F}$,\\
(ii) $\mathcal{C}$ does not contain any line and
the image of any color of $\mathcal{F}$ is zero.
\end{defi}
One can define a {\it colored fan} as a set of colored cones such that any two of them intersect in a common colored face (see \cite[def.1.14]{Pa06} for the precise definition).

Then $G/H$-embeddings are classified in terms of colored fans. 
Define a {\it simple} embedding as an open $G$-stable subvariety of $X$ containing exactly one closed $G$-orbit.  Let $X$ be a $G/H$-embedding and $\mathbb{F}$ its colored fan. Then $X$ is covered by its simple subembeddings, and each of them corresponds to a maximal colored cone of $\mathbb{F}$. 
Denote by $\mathcal{D}_X$ the set of simple roots in $S\backslash I$ which correspond to colors of $\mathbb{F}$.

\section{Classification of smooth projective embeddings with Picard number~1}\label{2}

The Picard number $\rho$ of a smooth $G/H$-embedding $X$ satisfies 
$$\rho=r+\sharp(S\backslash I)-\sharp(\mathcal{D}_X)$$
where $r$ is the number of rays of the colored fan of $X$ minus the rank $n$ \cite[(4.5.1)]{Pa06}. Since $X$ is projective, its colored fan is complete ({\it i.e.} it covers $N_\Rbb$) and then $r\geq 1$. Moreover $\mathcal{D}_X\subset S\backslash I$, so $\rho=1$ if and only if $r=1$ and $\mathcal{D}_X=S\backslash I$. In particular the colored fan of $X$ has exactly $n+1$ rays.

\begin{lem}\label{unique}
Let $G/H$ be a horospherical homogeneous space.  Up to isomorphism of varieties, there exists at most one smooth projective $G/H$-embedding with Picard number~1.
\end{lem}
\begin{proof}
Let $X$ and $X'$ be two smooth projective $G/H$-embeddings with respective colored fans $\mathbb{F}$ and $\mathbb{F'}$ and with Picard number~1.  Denote by $e_1,\ldots,e_{n+1}$ the primitive elements of the $n+1$ rays of $\mathbb{F}$. By the smoothness criterion of \cite[Chap.2]{Pa06}, $(e_1,\ldots,e_n)$ is a basis of $N$, $e_{n+1}=-e_1-\cdots -e_n$, the images in $N$ of the colors are disjoint and contained in  $\{e_1,\ldots,e_{n+1}\}$. The same happens for $\mathbb{F'}$.  Then there exists an automorphism $\phi$ of the lattice $N$ which stabilizes each image $\check\alpha_M$ of a color and satisfies $\mathbb{F}=\phi(\mathbb{F'})$. Thus the varieties $X$ and $X'$ are isomorphic  \cite[Prop. 3.10]{Pa06}.
\end{proof}

 If it exists, we call $X^1$ the unique smooth projective $G/H$-embedding with Picard number~1 and we say that $G/H$ is "special".

\subsection{Projective space}

The following result asserts in particular that $X^1$ is a projective space when $n\geq 2$.

\begin{teo}\label{thesproj}
Let $G/H$ be a "special" homogeneous space.
Then $X^1$ is isomorphic to a projective space in the following cases:
\begin{description}
\item (i) $\sharp(\mathcal{D}_{X^1})\leq n$,
\item (ii) $n\geq 2$
\item (iii) $n=1$, $\sharp(\mathcal{D}_{X^1})=2$ and the two simple roots of $\mathcal{D}_{X^1}$ are not in the same connected component of $\Gamma$.
\end{description}
\end{teo}

\begin{proof}
(i) In that case, a maximal colored cone of the colored fan of $X$ contains all colors. Then the corresponding simple subembedding of $X^1$, whose closed orbit is a point \cite[Lem.2.8]{Pa06}, is affine \cite[th.3.1]{Kn91} and smooth. So it is necessarily a horospherical $G$-module $V$ \cite[Lem.2.10]{Pa06}. Thus $\Pbb(\Cbb\oplus V)$ is a smooth projective $G/H$-embedding with Picard number~1. Then by Lemma \ref{unique}, $X^1$ is isomorphic to $\Pbb(\Cbb\oplus V)$.\\

(ii) We may assume that $\sharp(\mathcal{D}_{X^1})=n+1$.
Denote by by $\alpha_1,\ldots,\alpha_{n+1}$ the elements of $S\backslash I$ and by $\Gamma_i$ the Dynkin diagram $\Gamma_{S\backslash\{\alpha_i\}}$. The smoothness criterion of horospherical varieties \cite[Chap.2]{Pa06} applied to $X^1$ tells us two things.

First, for all $i\in\{1,\ldots,n+1\}$ and for all $j\neq i$, $\alpha_j$ is a simple end ("simple" means not adjacent to a double edge) of a connected component $\Gamma_i^j$ of $\Gamma_i$ of type $A_m$ or $C_m$. Moreover the $\Gamma_i^j$ are distinct, in other words, each connected component of $\Gamma_i$ has at most one vertex among the $(\alpha_i)_{i\in\{1,\ldots,n+1\}}$.

Secondly, $(\check\alpha_{iM})_{i\in\{1,\ldots,n\}}$ is a basis of $N$ and $\check\alpha_{(n+1)M}=-\check\alpha_{1M}-\cdots -\check\alpha_{nM}$. Thus a basis of  $M$ (dual of $N$) is of the form $$(\omega_{\alpha_i}-\omega_{\alpha_{n+1}}+\chi_i)_{i\in\{1,\ldots,n\}}$$ where $\chi_i$ is a character of the center $C$ of $G$,  for all $i\in\{1,\ldots,n\}$.

Let us prove that a connected component of $\Gamma$ contains at most one vertex among the $(\alpha_i)_{i\in\{1,\ldots,n+1\}}$. Suppose the contrary:  there exist $i,j\in\{1,\ldots,n+1\}$, $i\neq j$ such that $\alpha_i$ and $\alpha_j$ are vertices of a connected component of $\Gamma$. One can choose $i$ and $j$ such that there is no vertex among the $(\alpha_k)_{k\in\{1,\ldots,n+1\}}$ between $\alpha_i$ and $\alpha_j$. Since $n\geq 2$, there exists an integer $k\in\{1,\ldots,n+1\}$ different from $i$ and $j$. Then we observe that $\Gamma_k$ does not satisfy the condition that each of its connected component has at most one vertex among the $(\alpha_i)_{i\in\{1,\ldots,n+1\}}$ (because $\Gamma_k^i=\Gamma_k^j$).

Thus we have proved that  
\begin{equation}\label{eq1}
\Gamma=\bigsqcup_{j=0}^{n+1}\Gamma^j
\end{equation} such that for all $j\in\{1,\ldots,n+1\}$, $\Gamma^j$ is a connected component of $\Gamma$ of type $A_m$ or $C_m$ in which $\alpha_j$ is a simple end.

For all $\lambda\in\Lambda^+$, denote by by $V(\lambda)$ the simple $G$-module of weight $\lambda$. Then Equation \ref{eq1} tells us that the projective space $$\Pbb(V(\omega_{\alpha_{n+1}})\oplus V(\omega_{\alpha_1}+\chi_1)\oplus\cdots\oplus V(\omega_{\alpha_n}+\chi_n))$$  is a smooth projective $G/H$-embedding with Picard number~1. Thus $X^1$ is isomorphic to this projective space.\\

(iii)  As in case (ii), one checks that $X^1$ is isomorphic to $\Pbb(V(\omega_{\alpha_{2}})\oplus V(\omega_{\alpha_1}+\chi_1))$ (where $\chi_1$ is a character of $C$).
\end{proof}

\subsection{When $X^1$ is not isomorphic to a projective space}\label{2.2}

According to Theorem \ref{thesproj} we have to consider the case where the rank of $G/H$ is~1 and where there are two colors corresponding to  simple roots $\alpha$ and $\beta$  in the same connected component of $\Gamma$.
As we have seen in the proof of Theorem \ref{classif}, the lattice $M$ (here of rank~1) is generated by the character $\omega_\alpha-\omega_\beta+\chi$ where $\chi$ is a character of the center $C$ of $G$. Moreover, $H$ is the kernel of $\omega_\alpha-\omega_\beta+\chi:\, P(\omega_\alpha)\cap P(\omega_\beta)\lra\Cbb^*$. 

Let us reduce to the case where $G$ is semi-simple.

\begin{prop}
Let  $H'=G'\cap H$. Then $G/H$ is isomorphic to $G'/H'$. 
\end{prop}

\begin{proof}
We are going to prove that $G/H$ and $G'/H'$ are both isomorphic to a horospherical homogeneous space $(G'\times \Cbb^*)/H''$.
In fact $G/H$ is isomorphic to $(G'\times P/H)/\tilde{H}$ \cite[Proof of Prop.3.10]{Pa06}, where $$\tilde{H}=\{(g,pH)\in G'\times P/H\mid gp\in H\}.$$ Similarily $G'/H'$ is isomorphic to $(G'\times P'/H')/\tilde{H'}$ where $P'=P\cap G'$ and $\tilde{H'}$ defined as the same way as $\tilde{H}$.
Moreover the morphisms $$\begin{array}{ccccccc}
P/H & \lra & \Cbb^*  & \mbox{ and } &  P'/H' & \lra & \Cbb^*\\
pH & \lmt & (\omega_\alpha-\omega_\beta+\chi)(p) &\mbox{~~~~~} & p'H' & \lmt & (\omega_\alpha-\omega_\beta)(p')
\end{array}$$
are isomorphisms. Then $$\displaylines{\tilde{H}=\{(p',c)\in P'\times\Cbb^*\mid (\omega_\alpha-\omega_\beta+\chi)(p')=c^{-1}\}\hfill\cr\hfill =\{(p',c)\in P'\times\Cbb^*\mid (\omega_\alpha-\omega_\beta)(p')=c^{-1}\}=\tilde{H'}.}$$ 
This completes the proof.
\end{proof}

\begin{rem}\label{remP/H}
In fact $P/H\simeq\Cbb^*$ acts on $G/H$ by right multiplication, so it acts on the $\Cbb^*$-bundle $G/H\lra G/P$ by homotheties on fibers.  Moreover, this action extends to $X^1$.
\end{rem}

We may assume that $G$ is semi-simple. Let $G_1,\ldots,G_k$ the normal simple subgroups of $G$, so that $G$ is the quotient of the product $G_1\times\cdots\times G_k$ by a central finite group $C_0$. We can suppose that $C_0$ is trivial because $G/H\simeq G.C_0/H.C_0$. If $\alpha$ and $\beta$ are simple roots of the connected component corresponding to $G_i$, denote by $H_i$ is the kernel of $\omega_\alpha -\omega_\beta$ in the parabolic subgroup $P(\omega_\alpha)\cap P(\omega_\beta)$ of $G_i$. Then $$H=G_1\times\cdots\times G_{i-1}\times H_i\times G_{i+1}\times\cdots\times G_k$$ and  $G/H=G_i/H_i$.

So from now, without loss of generality, we suppose that $G$ is simple.

\begin{teo}\label{classif}
With the assumptions above,  $G/H$ is "special" if and only if $(\Gamma,\alpha,\beta)$ appears in  the following list (up to exchanging $\alpha$ and $\beta$).\footnote{The numerotation of the simple roots is that of \cite{Bo75}.}
\begin{enumerate}
\item $(A_m,\alpha_1,\alpha_m)$

\item $(A_m,\alpha_i,\alpha_{i+1})$ with $i\in\{1,\ldots,m-1\}$

\item $(B_m,\alpha_{m-1},\alpha_m)$

\item $(B_3,\alpha_1,\alpha_3)$

\item $(C_m,\alpha_i,\alpha_{i+1})$ with $i\in\{1,\ldots,m-1\}$

\item $(D_m,\alpha_{m-1},\alpha_m)$

\item $(F_4,\alpha_2,\alpha_3)$

\item $(G_2,\alpha_2,\alpha_1)$
\end{enumerate}
\end{teo}

\begin{proof}
The Dynkin diagrams $\Gamma_{S\backslash\{\alpha\}}$ and $\Gamma_{S\backslash\{\beta\}}$ respectively. They are of type $A_m$ or $C_m$ by the smoothness criterion \cite[Chap.2]{Pa06}. Moreover $\alpha$ and $\beta$ are simple ends of $\Gamma_{S\backslash\{\beta\}}$ and $\Gamma_{S\backslash\{\alpha\}}$ respectively.

Suppose $\Gamma$ is of type $A_m$. If $\alpha$  equals $\alpha_1$ then, looking at $\Gamma_{S\backslash\{\alpha\}}$, we remark that $\beta$ must be $\alpha_2$ or $\alpha_m$. So we are in Case 2 or 1. If $\alpha$ equals $\alpha_m$ the argument is similar. Now if $\alpha$ is not an end of $\Gamma$, in other words if $\alpha=\alpha_i$ for some $i\in\{2,\ldots,m-1\}$ then, looking at $\Gamma_{S\backslash\{\alpha\}}$, we see that $\beta$ can be $\alpha_1$, $\alpha_{i-1}$, $\alpha_{i+1}$ or $\alpha_m$. The cases where $\beta$ equals $\alpha_1$ or $\alpha_m$ are already done and the case where $\beta$ equals $\alpha_{i-1}$ or $\alpha_{i+1}$ is Case 2.

The study of the other cases is analogous and left to the reader.
\end{proof}

In the next two sections we are going to study the variety $X^1$ for each case of this theorem. In particular we will see that $X^1$ is not isomorphic to a projective space in all cases.

\section{Homogeneous varieties}\label{3}

In this section, with the notation of Section \ref{2.2}, we are going to prove that $X^1$ is homogeneous in Cases 1, 2 and 6.\\

In all cases (1 to 8), there are exactly 4 projective $G/H$-embeddings and there are all smooth; they correspond to the 4 colored fans consisting of the two half-lines of $\Rbb$, without color, with one of the two colors and with the two colors, respectively (see \cite[Ex.1.19]{Pa06} for a similar example).

Let us realize $X^1$ in a projective space as follows.  The homogeneous space $G/H$ is isomorphic to the orbit of the point $[v_{\omega_\beta}+v_{\omega_\alpha}]$ in  $\Pbb(V(\omega_\beta)\oplus V(\omega_\alpha))$, where $v_{\omega_\alpha}$ and $v_{\omega_\beta}$ are highest weight vectors of $V(\omega_\alpha)$ and $V(\omega_\beta)$ respectively. Then $X^1$ is the closure of this orbit in $\Pbb(V(\omega_\beta)\oplus V(\omega_\alpha))$, because both have the same colored cone ({\it i.e.} that with two colors)\footnote{See \cite[Chap.1]{Pa06} for the construction of the colored fan of a $G/H$- embbedding.}. 

We will describe the other $G/H$-embeddings later in the proof of Lemma \ref{lemaut}.

\begin{prop}
In Case 1, $X^1$ is isomorphic to the quadric $Q^{2m}=\operatorname{SO}_{2m+2}/P(\omega_{\alpha_1})$.  
\end{prop}

\begin{proof}
Here, the fundamental representations $V(\omega_\alpha)$ and $V(\omega_\beta)$ are the simple $\operatorname{SL}_{m+1}$-modules $\Cbb^{m+1}$ and its dual $(\Cbb^{m+1})^*$, respectively. Let denote by $Q$ the quadratic form on $\Cbb^{m+1}\oplus(\Cbb^{m+1})^*$ define by $Q(u,u^*)=\langle u^*,u\rangle$. Then $Q$ is invariant under the action of $\operatorname{SL_{m+1}}$. Moreover $Q(v_{\omega_\alpha}+v_{\omega_\beta})=0$, so that $X^1$ is a subvariety of the quadric ($Q=0$) in $\Pbb(\Cbb^{m+1}\oplus(\Cbb^{m+1})^*)=\Pbb(\Cbb^{2m+2})$.

We complete the proof by computing the dimension of $X^1$: \begin{equation}\label{dim}
\operatorname{dim}\,X^1=\operatorname{dim}\,G/H=1+\operatorname{dim}\,G/P=1+\sharp(R^+\backslash R_I^+)
\end{equation}
where $R^+$ is the set of positive roots of $(G,B)$ and $R_I^+$ is the set of positive roots generated by simple roots of $I$.
So $\operatorname{dim}\,X^1=\operatorname{dim}\,Q^{2m}=2m$ and $X^1=Q^{2m}$.
\end{proof}

\begin{prop}\label{cas2}
In Case 2, $X^1$ is isomorphic to the grassmannian $Gr_{i+1,m+2}$.
\end{prop}

\begin{proof}
The fundamental representations of $\operatorname{SL_{m+1}}$ are the simple $\operatorname{SL_{m+1}}$-modules
$$\begin{array}{c} V(\omega_{\alpha_i})=\bigwedge^i\Cbb^{m+1}
\end{array}$$
 and a highest weight vector of $V(\omega_{\alpha_i})$ is $e_1\wedge\cdots\wedge e_i$ where $e_1,\ldots,e_{m+1}$ is a basis of $\Cbb^{m+1}$.

 We have
$$\xymatrix{
    X\, \ar@{^{(}->}[r]  & \,\Pbb(\bigwedge^i\Cbb^{m+1}\oplus\bigwedge^{i+1}\Cbb^{m+1})\\
    G/H\, \ar@{=}[r] \ar@{^{(}->}[u]  & G.[e_1\wedge\cdots\wedge e_i+e_1\wedge\cdots\wedge e_{i+1}] \ar@{^{(}->}[u]
  }$$
Complete $(e_1,\ldots,e_{m+1})$ to obtain a basis $(e_0,\ldots,e_{m+1})$ of $\Cbb^{m+2}$, then the morphism
$$\begin{array}{ccc}
\bigwedge^i\Cbb^{m+1}\oplus\bigwedge^{i+1}\Cbb^{m+1} & \lra & \bigwedge^{i+1}\Cbb^{m+2}\\
x+y & \lmt & x\wedge e_0 +y
\end{array}$$
 is an isomorphism.
Then $X^1$ is a subariety of the grassmannian $$\begin{array}{c} G_{i+1,m+2}\simeq \operatorname{SL}_{m+2}.[e_1\wedge\cdots\wedge e_i\wedge(e_0+e_{i+1})]\subset\Pbb(\bigwedge^{i+1}\Cbb^{m+1}).\end{array}$$ 
Moreover the dimension of $X^1$ is the same as the dimension of $G_{i+1,m+2}$. Indeed, one checks that  $\operatorname{dim}\,X^1=(i+1)(m+1-i)$ using Formula \ref{dim}.
\end{proof}

\begin{prop}
In Case 6, $X^1$ is isomorphic to the spinor variety $\operatorname{SO}(2m+1)/P(\omega_{\alpha_m})$.
\end{prop}

\begin{proof}
The direct sum $V(\omega_\alpha)\oplus V(\omega_\beta)$ of the two half-spin representations of $\operatorname{SO}(2m)$ is isomorphic to the spin representation of $\operatorname{SO}(2m+1)$.  Moreover $v_{\omega_\alpha}+v_{\omega_\beta}$ is in the orbit of a highest weight vector of the spin representation of $\operatorname{SO}(2m+1)$. Thus we deduce that $X^1$ is a subvariety of $\operatorname{SO}(2m+1)/P(\omega_{\alpha_m})$.  Then we complete the proof by comparing dimensions. In fact $\operatorname{dim}\,X^1=\frac{m(m+1)}{2}$ (by Formula \ref{dim}).
\end{proof}

\section{Non-homogeneous varieties}\label{4}

With the notation of Section \ref{2.2} we prove in this section the following result.

\begin{teo}\label{mainteo}
In Cases 3, 4, 5, 7 and 8 (and only in these cases), $X^1$ is not homogeneous. 

Moreover the automorphism group of $X^1$ is $(\operatorname{SO}(2m+1)\times\Cbb^*)\ltimes V(\omega_{\alpha_m})$, $(\operatorname{SO}(7)\times\Cbb^*)\ltimes V(\omega_{\alpha_3})$, $((\operatorname{Sp}(2m)\times\Cbb^*)/\{\pm 1\})\ltimes V(\omega_{\alpha_1})$, $(\operatorname{F}_4\times\Cbb^*)\ltimes V(\omega_{\alpha_4})$ and $(\operatorname{G}_2\times\Cbb^*)\ltimes V(\omega_{\alpha_1})$ respectively.

Finally, $X^1$ has two orbits under its automorphism group.
\end{teo}

One can remark that in Case 5, Theorem \ref{mainteo} follows from \cite[Chap.3 and Prop.5.1]{Mi05} and the following result.

\begin{prop}\label{oddgrass}
In Case 5, $X^1$ is isomorphic to the odd symplectic grassmannian $Gr_{i+1,2m+1}$.
\end{prop}

\begin{proof}
As in the proof of Proposition \ref{cas2},  $X^1$ is a subvariety of the odd symplectic grassmannian $$G_{i+1,2m+1}\simeq\overline{\operatorname{Sp}_{2m+1}.[e_1\wedge\cdots\wedge e_i\wedge(e_0+e_{i+1})]}\subset\Pbb(\Lambda^{i+1}\Cbb^{m+1}).$$
Again we complete the proof by comparing the dimensions. Indeed, by Formula \ref{dim} $$\operatorname{dim}\,X^1=(i+1)(2m-i)-\frac{i(i+1)}{2}=\operatorname{dim}\,Gr_{i+1,2m+1}\,\,\, \cite[\mbox{Prop 4.1}]{Mi05}.$$
\end{proof}

Now, let $X$ be one of the varieties $X^1$ in Cases 3, 4, 7 and 8.

Then $X$ has three orbits under the action of $G$ (the open orbit $X_0$ isomorphic to $G/H$ and two closed orbits). Recall that $X$ can be seen as a subvariety of $\Pbb(V(\omega_\alpha)\oplus V(\omega_\beta))$.  Let $P_Y:=P(\omega_\alpha)$, $P_Z:=P(\omega_\beta)$ and denote by $Y$ and $Z$ the closed orbits, isomorphic to $G/P_Y$ and $G/P_Z$ respectively. (In Case 8 where $G$ is of type $G_2$, we have $\alpha=\alpha_2$ and $\beta=\alpha_1$.)

Let $X_Y$ be  the simple $G/H$-embedding of $X$, we have $X_Y=X_0\cup Y$. Then $X_Y$ is a homogeneous vector bundle over $G/P_Y$ \cite[Chap.2]{Pa06} in the sense of the following definition.
\begin{defi}
Let $P$ be a parabolic subgroup of $G$ and $V$ a $P$-module. Then the homogeneous vector bundle $G\times^PV$ over $G/P$ is the quotient of the direct product $G\times V$ by the equivalence relation $\sim$ defined by $$\forall g\in G,\,\forall p\in P, \forall v\in V,\quad (g,v)\sim(gp^{-1},p.v).$$
\end{defi}
Specifically, $X_Y=G\times^{P_Y}V_Y$ where $V_Y$ is a simple $P_Y$-module of highest weight $\omega_\beta-\omega_\alpha$, and similarily, $X_Z=G\times^{P_Z}V_Z$ where $V_Z$ is a simple $P_Z$-module of highest weight $\omega_\alpha-\omega_\beta$.

Denote by $\Aut(X)$ the automorphism group of $X$ and  $\Aut ^0 (X)$ the connected component of $\Aut (X)$ containing the identity.  
\begin{rem}
Observe that $\Aut(X)$ is a linear algebraic group. Indeed $\Aut(X)$ acts on the Picard group of $X$ which equals $\Zbb$ (the Picard group of a horospherical variety is free \cite{Br89}). This action is necessarily trivial. Then $\Aut(X)$ acts on the projectivization of the space of global sections of a very ample bundle. This gives a faithful representation of $\Aut(X)$.
\end{rem}

We now complete the proof of Theorem \ref{mainteo} by proving several lemmas.

\begin{lem}
The closed orbit $Z$ of $X$ is stable under $\Aut ^0 (X)$.
\end{lem}

\begin{proof}
We are going to prove that the normal bundle $N_Z$ of $Z$ in $X$ has no nonzero global section. This will imply that $\Aut^0(X)$ stabilises $Z$, because the Lie algebra $\Lie(\Aut^0(X))$ is the space of global sections $H^0(X,T_{X})$ of the tangent bundle $T_X$ of $X$ \cite[Chap.2.3]{Ak95} and we have the following exact sequence
$$0\lra T_{X,Z}\lra T_{X}\lra N_Z\lra 0$$
where $T_{X,Z}$ is the subsheaf of $T_X$ consisting of vector fields that vanish along $Z$. Moreover $H^0(X,T_{X,Z})$ is the Lie algebra of the subgroup of $\Aut^0(X)$ that stabilizes $Z$.

The total space of $N_Z$ is the vector bundle $X_Z$. So using the Borel-Weil theorem \cite[4.3]{Ak95}, $H^0(G/P_Z,N_Z)=0$ if and only if the smallest weight of $V_Z$ is not antidominant. The smallest weight of $V_Z$ is  $w_0^\beta(\omega_\alpha-\omega_\beta)$ where $w_0^\beta$ is the longest element of $W_{S\backslash\{\beta\}}$.  
Let $\gamma\in S$, then $$\langle w_0^\beta(\omega_\alpha-\omega_\beta),\check\gamma\rangle=\langle\omega_\alpha-\omega_\beta,w_0^\beta(\check\gamma)\rangle.$$
If $\gamma$ is different from $\beta$  then $w_0^\beta(\check\gamma)=-\check\delta$ for some $\delta\in S\backslash\{\beta\}$. So we only have to compute $w_0^\beta(\check\beta)$. 

In Case 3, $\beta=\alpha_m$, so $w_0^\beta$ maps $\alpha_i$ to $-\alpha_{m-i}$ for all $i\in\{1,\ldots,m-1\}$\footnote{The numerotation of simple roots is still that of \cite{Bo75}.}. Here, $\omega_\beta=\alpha_1+2\alpha_2+\cdots+m\alpha_m$. Then, using the fact that $w_0^\beta(\omega_\beta)=\omega_\beta$, we have  $w_0^\beta(\beta)=\alpha_1+\cdots+\alpha_m$ so that $w_0^\beta(\check\beta)=2(\alpha_1+\cdots+\alpha_{m-1})+\alpha_m$ and $\langle\omega_\alpha-\omega_\beta,w_0^\beta(\check\gamma)\rangle=1$ (because $\alpha=\alpha_{m-1}$).

 The computation of $w_0^\beta(\check\beta)$ in the other cases is similar and left to the reader. In all four cases,$\langle\omega_\alpha-\omega_\beta,w_0^\beta(\check\beta)\rangle>0$ (this equals 1 in Cases 3, 4, 7 and 2 in Case 8).
This proves that $w_0^\beta(\omega_\alpha-\omega_\beta)$ is not antidominant.
\end{proof}

\begin{rem}\label{sectionNY}
Using the Borel-Weil theorem we can also compute $H^0(G/P_Y,N_Y)$ by the same method. We find that in Cases 3, 4, 7 and 8, this $G$-module is isomorphic to the simple $G$-module $V(\omega_\beta)$, $V(\omega_\beta)$, $V(\omega_{\alpha_4})$ and $V(\omega_\beta)$ respectively.   
\end{rem}

We now prove the following lemma.

\begin{lem}\label{lemaut}
In Cases 3, 4, 7 and 8, $\Aut^0(X)$ is $(\operatorname{SO}(2m+1)\times\Cbb^*)\ltimes V(\omega_{\alpha_m})$, $(\operatorname{SO}(7)\times\Cbb^*)\ltimes V(\omega_{\alpha_3})$, $(\operatorname{F}_4\times\Cbb^*)\ltimes V(\omega_{\alpha_4})$ and $(\operatorname{G}_2\times\Cbb^*)\ltimes  V(\omega_{\alpha_1})$ respectively.
\end{lem}
\begin{rem}
By Remark \ref{remP/H}, we already know that the action of $G$ on $X$ extends to an action of $G\times\Cbb^*$($\simeq G\times P/H$). Moreover $\tilde{C}:=\{(c,c^{-1}H)\in C\times P/H\}\subset G\times\Cbb^*$ acts trivially on $X$.
\end{rem}
\begin{proof}

Let $\pi:\tilde{X}\lra X$ be the blowing-up of $Z$ in $X$.  Since $Z$ and $X$ are smooth, $\tilde{X}$ is smooth; it is also a projective $G/H$-embbeding. In fact $\tilde{X}$ is the projective bundle $$\phi:\,G\times^{P_Y}\Pbb(V_Y\oplus\Cbb)\lra G/P_Y$$ where $P_Y$ acts trivially on $\Cbb$ (because both are projective $G/H$-embeddings with exactly the same color). Moreover the exceptional divisor $\tilde{Z}$ of $\tilde{X}$ is $G/P$.

Let us remark that $\Aut^0(\tilde{X})$ is isomorphic to $\Aut^0(X)$. Indeed, it contains $\Aut^0(X)$ because $Z$ is stable under the action of $\Aut^0(X)$ and we know, by a result of A.~Blanchard, that $\Aut^0(\tilde{X})$ acts on $X$ such that $\pi$ is equivariant \cite[Chap.2.4]{Ak95}. 

Now we are going to compute $\Aut^0(\tilde{X})$. 
Observe that $H^0(G/P_Y,N_Y)$ acts on $\tilde{X}$ by translations on the fibers of $\phi$:
$$\forall s\in H^0(G/P_Y,N_Y),\, \forall (g_0,[v_0,\xi])\in G\times^{P_Y}\Pbb(V_Y\oplus\Cbb),\,
s.(g_0,[v_0,\xi])=(g_0,[v_0+\xi v(g_0),\xi])$$ where $v(g_0)$ is the element of $V_Y$ such that $s(g_0)=(g_0,v(g_0))$.

Then the group $((G\times\Cbb^*)/\tilde{C})\ltimes H^0(G/P_Y,N_Y)$ acts effectively on $\tilde{X}$
(the semi-product is defined by $((g',c'),s').((g,c),s)=((g'g,c'c),c'g's+s')$). 
In fact we are going to prove that $$\Aut^0(\tilde{X})=((G\times\Cbb^*)/\tilde{C})\ltimes H^0(G/P_Y,N_Y).$$
 
By \cite[Chap.2.4]{Ak95} we know that $\Aut^0(\tilde{X})$ exchanges the fibers of $\phi$ and induces an automorphism of $G/P_Y$. Moreover we have $\Aut^0(G/P_Y)=G/C$  in our four cases \cite[Chap.3.3]{Ak95}.
So we have an exact sequence
$$0\lra A\lra \Aut^0(\tilde{X})\lra G/C\lra 0$$
where $A$ is the set of automorphisms which stabilise each  fiber of the projective bundle $\tilde{X}$. In fact $A$ consists of affine transformations in fibers.

Then $H^0(G/P_Y,N_Y)$ is the subset of $A$ consisting of translations. Let $A_0$ be the subset of $A$ consisting of linear transformations in fibers. Then $A_0$ fixes $Y$ so that $A_0$ acts on the blowing-up  $\tilde{\tilde{X}}$ of $Y$ in $\tilde{X}$. Moreover $\tilde{\tilde{X}}$ is a $\Pbb^1$-bundle over $G/P$, it is in fact the toroidal $G/H$-embedding ({\it i.e.} without colors) \cite[Ex.1.13 (2)]{Pa06}. As before we know that $\Aut^0(\tilde{\tilde{X}})$ exchanges the fibers of that $\Pbb^1$-bundle and induces an automorphism of $G/P$. Moreover we have $\Aut^0(G/P)=G/C$ and then $\Aut^0(\tilde{\tilde{X}})=(G\times\Cbb^*)/\tilde{C}$. We deduce that $A_0=\Cbb^*$.\\

We complete the proof by Remark \ref{sectionNY}.

\end{proof}
\begin{rem}
We can use the same arguments for the odd symplectic grassmannian (Case 5). Then $Y$ is stabilised by $\Aut^0(X)$ and $H^0(G/P_Z,N_Z)=V(\omega_{\alpha_1})\simeq\Cbb^{2m}$. We deduce that $\Aut^0(X)=((\operatorname{Sp}(2m)\times\Cbb^*)/\{\pm 1\})\ltimes\Cbb^{2m}$.
\end{rem}
 
To complete the proof of Theorem \ref{mainteo} we prove the following lemma.

\begin{lem}
The automorphism group of $X$ is connected. 
\end{lem}

\begin{proof}
Let $\phi$ be an automorphism of $X$. We want to prove that $\phi$ is in $\Aut^0(X)$.  It acts by conjugation on $\Aut^0(X)$. Let $L$ be a Levi subgroup of $\Aut^0(X)$. Then $\phi^{-1}L\phi$ is again a Levi subgroup of $\Aut^0(X)$. But all Levi subgroups are conjugated in $\Aut^0(X)$. So we can suppose, without loss of generality, that $\phi$ stabilises $L$. 

Then $\phi$ induces an automorphism of the direct product of $\Cbb^*$ with a simple group $G$ of type $B_m$, $C_m$, $F_4$ or $G_2$.  It also induces an automorphism of $G$ which is necessarily an inner automorphism of $G$ (because there is no non-trivial automorphism of the Dynkin diagram of $G$). So we can assume now that $\phi$ commutes with all elements of $G$.

Then $\phi$ stabilises the open orbit $G/H$ of $X$. Let $x_0:=H\in G/H\subset X$ and $x_1:=\phi(x_0)\in G/H$. Since $\phi$ commutes with the elements of $G$, the stabilizer of $x_1$ is also $H$. So $\phi$ acts on $G/H$ as an element of $N_G(H)/H=P/H\simeq \Cbb^*$ (where $N_G(H)$ is the normalizer of $H$ in $G$).  Then $\phi$ is an element of $\Cbb^*\subset\Aut^0(X)$.
\end{proof}

\section{On two-orbit varieties}\label{5}

Here we give a characterization of the two-orbit varieties obtained in Section \ref{4}.

\begin{teo}\label{caracterisation}
Let $X$ be a smooth projective variety with Picard number~1 and put $\mathbf{G}:=\Aut^0(X)$. Suppose that $\mathbf{G}$ is not semi-simple and that $X$ has two orbits under the action of $\mathbf{G}$. Denote by $Z$ the closed orbit. 

Then the codimension of $Z$ is at least $2$.

Suppose in addition that the blowing-up of $X$ along $Z$ has also two orbits under the action of $\mathbf{G}$. Then $X$ is one of the varieties $X^1$ obtained in the cases 3, 4, 5, 7 and 8.
\end{teo}

\begin{rem}
The converse implication holds by Section \ref{4}.
\end{rem}

To prove this result we need a result of D.~Akhiezer on two-orbits varieties with codimension one closed orbit.

\begin{lem}[Th.1 of \cite{Ah83}]\label{Akhiezer}

Let $X$ a smooth complete variety with an effective action of the (connected linear non semi-simple) group $\mathbf{G}$.
Suppose that $X$ has two orbits under the action of $\mathbf{G}$ and that the closed orbit is of codimension $1$.
Let $G$ be a maximal semi-simple subgroup of $\mathbf{G}$.\\

Then there exist a parabolic subgroup $P$ of $G$ and a $P$-module $V$ such that:\\
(i) the action of $P$ on $\Pbb(V)$ is transitive;\\
(ii) there exists an irreducible $G$-module $W$ and a surjective $P$-equivariant morphism $W\lra V$;\\
(iii) $X=G\times^P\Pbb(V\oplus\Cbb)$.\\

In particular, $X$ is a horospherical $G$-variety of rank $1$.
\end{lem}

\begin{rems}
It follows from (i) that $V$ is an (irreducible) horospherical $L$-module of rank $1$, where $L$ is a Levi subgroup of $G$. This is why $X$ is horospherical of rank $1$.

The $G$-module $W$ is the set of global sections of the vector bundle $G\times^P V$.
\end{rems}
 
\begin{proof}[Proof of Theorem \ref{caracterisation}] 
If $Z$ is of codimension $1$,  Lemma \ref{Akhiezer} tells us that $X$ is a horospherical $G$-variety. Moreover the existence of the irreducible $G$-stable divisor $Z$ tells us that one of the two rays of the colored fan of $X$ has no color \cite[Chap 1]{Pa06}. Then $X$ satisfies the condition (i) of Theorem \ref{thesproj} and $X=\Pbb(V\oplus\Cbb)$.
Since $X$ is not homogeneous we conclude that the codimension of $Z$ is at least~2.

Denote by $\tilde{X}$ the blowing-up of $X$ along $Z$. Then $\tilde{X}$ is a horospherical $G$-variety by Lemma \ref{Akhiezer}. Moreover $X$ and $\tilde{X}$ have the same open $G$-orbit, so that $X$ is also a horospherical $G$-variety.
We conclude by Theorem \ref{mainteo}.
\end{proof}

\end{document}